\documentclass[10pt]{article}

\usepackage[T1]{fontenc}
\usepackage{array,amsfonts}
\usepackage{amstext}
\usepackage{amsmath,amsthm}
\usepackage{amssymb}
\usepackage{graphicx}
\usepackage{textcomp}
\usepackage{graphics,xcolor}
\usepackage{stmaryrd}

\newcommand\ds{\displaystyle}
\newcommand\tn{\textnormal}

\newtheorem{theo}{Theorem}

\newtheorem{lemme}{Lemma}
\newtheorem{rem}{Remark}

\title{Lipschitz stability in an inverse problem\\ for the wave equation}
\author{Lucie Baudouin$^{1,2,}$\footnote{e-mail: {\tt baudouin@laas.fr}}\\
{\it\footnotesize $^{1}$ CNRS ; LAAS ; 7 avenue du colonel Roche, F-31077 Toulouse Cedex 4, France}\\
{\it\footnotesize $^{2}$ Universit\'e de Toulouse ; UPS, INSA, INP, ISAE, UT1, UTM ; LAAS ; F-31077 Toulouse, France.}}
\date{November 2001}
\begin{document}

\maketitle

 \textbf{Abstract:} We are interested in the inverse problem of the determination  of the potential $p(x),~x\in\Omega\subset\mathbb{R}^n$ from the measurement of the normal derivative $\partial_\nu u$ on a suitable part $\Gamma_0$ of the boundary of $\Omega$, where
$u$ is the solution of the wave equation $\partial_{tt}u(x,t)-\Delta u(x,t)+p(x)u(x,t)=0$ set in $\Omega\times(0,T)$ and given Dirichlet boundary data.
More precisely, we will prove local uniqueness and stability for this inverse problem and the main tool will be a global Carleman estimate, result also interesting by itself.



\section{Introduction and main result}

Let $n \in \mathbb{N}$, $T>0$ and let $\Omega \subset \mathbb{R}^n$ be a
bounded domain with $C^2$-boundary $\partial\Omega$. Let $\Gamma_{0}$ be an
open subset of $\partial\Omega$. 
Throughout this paper, for a functional $v = v(x,t)$ with $x\in\Omega,~t\in (0,T)$, we use the following notations :
\begin{eqnarray*}
&&\nabla v = \left(\frac {\partial v}{\partial x_1},\dots, \frac {\partial v}{\partial x_n}\right),
\quad D^2v=\left(\frac {\partial^2 v}{\partial x_i\partial x_j}\right)_{1\leq i,j\leq n},
\quad \Delta v = \sum_{i=1}^{n}\frac {\partial^2 v}{\partial x_i^2},
\quad \partial_t v= \frac {\partial v}{\partial t},\\
&&\nu\in\mathbb{R}^n~\tn{ denotes the unit outward normal vector to } \partial\Omega~\tn{ and }
\partial_\nu v = \frac {\partial v}{\partial \nu}=\nabla v\cdot \nu.
\end{eqnarray*}

We consider the wave equation :
\begin{equation}\label{WE}
	\left\lbrace\begin{array}{ll}
		\partial_t^2y(x,t)-\Delta y(x,t)+q(x)y(x,t)=g(x,t),&\quad x\in\Omega,~t\in (0,T),\\
		y(x,t)=h(x,t),&\quad x \in\partial\Omega,~t\in (0,T),\\
		y(x,0)=y^0(x),\quad \partial_t y(x,0)= y^1(x),&\quad x\in\Omega.
	\end{array} \right.
\end{equation}
First of all, assuming that $y^0\in L^2(\Omega)$, $y^1\in H^{-1}(\Omega)$, $p\in L^\infty(\Omega)$, $h \in L^2(\partial\Omega \times (0,T))$  and $g \in L^1(0,T;L^2(\Omega))$ are known, and assuming the compatibility condition $h(x,0) = y_0(x)$ for all $x\in\partial\Omega$, the Cauchy problem \eqref{WE} is well-posed and one can also prove (using a method by transposition, since we have non-homogeneous boundary conditions) that
$$u\in C([0,T],L^2(\Omega)) \cap C^1([0,T],H^{-1}(\Omega)).$$
This result can be read in \cite{LionsMagenesBook} for instance. One will also find a classical existence and regularity result when the boundary data $h$ is equal to $0$ in Lemma~\ref{Energy} (see also \cite{Lions}), useful in the inverse problem result.\\

This paper treats at the same time two kinds of inverse problems which can be stated as follows.

\textbf{Non linear inverse Problem :}
Is it possible to retrieve the potential $q=q(x),~x\in\Omega~$ from measurement
of the  normal derivative $\left.\partial_\nu y\right|_{\Gamma_{0} \times (0,T)}$
where $y$ is the solution to (\ref{WE}), $\Gamma_0$ is a part large enough of the boundary $\partial\Omega$ and the observation time $T$ is also large enough ?\\

We will actually give local answer to this question. If we denote by $y[p]$ the weak solution of equation~\eqref{WE}, assuming that $p\in L^\infty(\Omega)$ is a given potential, we are concerned with the stability around $p$. That is to say $p$ and $y[p]$ are known while $q$ is unknown and we prove the following local lipschitz stability result.
In this direction, we will answer to two more precise problems.\\

\noindent\textbf{Uniqueness} : 
Under geometrical conditions on $\Gamma_{0}$ and $T$, does the equality 
$\partial_\nu y[q]=\partial_\nu  y[p]$ on $\Gamma_{0} \times(0,T)$
imply $q=p$ on $\Omega$ ?\\

\noindent\textbf{Stability} : 
Under geometrical conditions on $\Gamma_{0}$ and $T$, is it possible to estimate $\|q-p\|_{L^2(\Omega)}$ 
or better, a stronger norm of $(p-q)$, by a suitable norm of 
$\partial_\nu y[q]-\partial_\nu  y[p]$ on $\Gamma_{0}\times (0,T)$ ?\\

We will actually work on a linearized version of the inverse problem and consider the following wave equation :

\begin{equation}\label{WE2}
	\left\lbrace\begin{array}{ll}
		\partial_t^2u(x,t)-\Delta u(x,t)+q(x)u(x,t)=f(x)R(x,t),&\quad x\in\Omega,~t\in (0,T),\\
		u(x,t)=0,&\quad x \in\partial\Omega,~t\in (0,T),\\
		u(x,0)=0, \quad \partial_t u(x,0)= 0,&\quad x\in\Omega.
	\end{array} \right.
\end{equation}

\textbf{Linear inverse problem :}
Is it possible to determine $f(x),~x\in\Omega~$ from the knowledge of the normal derivative 
$\left.\partial_\nu u\right|_{\Gamma_{0} \times(0,T)} $ 
where $R$ and $q$ are given and $u$ is the solution to (\ref{WE2})?\\

These questions for the wave equation have all already received positive answers
since the uniqueness result for the linear inverse problem has been proved by M.V. \textsc{Klibanov} in \cite{Klibanov92} 
and Lipschitz stability results (for both linear and non-linear inverse problems) of M. \textsc{Yamamoto}, deriving from it, can be read in \cite{Yam99}.
The proof in \cite{Yam99} is based on a local Carleman estimate for the wave operator and a compactness-uniqueness argument in order to conclude to the stability from the uniqueness result and an observability estimate.
We aim in our document at giving a direct proof of a Lipschitz stability estimate from a global Carleman estimate, result also interesting by itself. Another gain of this new proof of the precise result given below is the weakened assumptions on the solution of the wave equation under study. Besides, from the Carleman estimate we prove in the sequel, we directly obtain that a measurement of the flux of the solution on a suitable part $\Gamma_0$ of the boundary  (instead of the whole boundary $\partial\Omega$) is sufficient.\\

To precisely state the results we will prove in this article, we introduce, for $m \geq 0$, the set
$$
	L^\infty_{\leq m} (\Omega) = \{q \in L^\infty(\Omega), \ s.t. \, \|q\|_{L^\infty(\Omega)} \leq m \}.
$$
Moreover, we also specify the geometrical assumption:
\begin{gather}
\exists x_{0}\not\in\Omega\quad \tn{such that} \quad \Gamma_{0}\supset \{x\in\partial\Omega ; (x-x_{0})\cdot \nu(x)\geq 0\}
\label{GCC}\\
T> \sup_{x\in \Omega} |x-x_0|
\label{GCCT}
\end{gather}

\begin{theo}\label{TWE}
Let $m>0$, $K>0$ and $r>0$. Let $p$ belong to  $L^\infty_{\leq m}(\Omega)$. 
Assume that the solution $y[p]$ of equation \eqref{WE} is such that 
\begin{equation}
\label{RegCond}
	\|y[p]\|_{H^1 (0,T;L^{\infty}(\Omega))} \leq K
\end{equation}
and assume also that the initial datum $y^0$ satisfies
\begin{equation}
	\label{InitialDataCond}
	 \inf \left\{|y^0(x)|,x\in (\Omega)\right\} \geq r.
\end{equation}
Under the hypothesis that $\Gamma_0$ satisfies the geometrical condition \eqref{GCC} and $T$ satisfies \eqref{GCCT},
then for all $q \in L^\infty_{\leq m}(\Omega)$, 
$ \partial_\nu y[q]-\partial_\nu y[p]$ belongs to $H^1(0,T;L^2(\Gamma_0)) $ 
and there exists a constant $C = C(m,T,K,r)>0$ such that for any $q \in L^\infty_{\leq m} (\Omega)$,
\begin{equation}
	\label{ReverseEst}
	\left\|\dfrac{\partial y[q]}{\partial\nu}-\dfrac{\partial y[p]}{\partial\nu}\right\|_{H^1(0,T;L^2(\Gamma_0))}
	\leq C \|p-q\|_{L^2(\Omega)},
\end{equation}
\begin{equation}
	\label{StabilityContinuous}
	\|q-p\|_{L^2(\Omega)}\leq C\left\|\dfrac{\partial y[q]}{\partial\nu}-\dfrac{\partial y[p]}{\partial\nu}\right\|_{H^1(0,T;L^2(\Gamma_0))} .
\end{equation}
\end{theo}
Estimate \eqref{StabilityContinuous} is the Lispchtiz stability of the inverse problem while \eqref{ReverseEst} gives the continuous dependance of the derivative of the solution with respect to the potential.
\begin{rem}\label{RemReg}
	The condition $y[p] \in H^1 (0,T;L^{\infty}(\Omega))$ can be guaranteed uniformly for $p \in L^\infty_{\leq m} (\Omega)$ with more constraints on the data $(y^0, y^1), \ g,\ h$ in \eqref{WE}. Indeed, if we assume $(y^0,y^1)  \in   H^2(\Omega) \times H^1(\Omega)$,
$g  \in  W^{1,1}(0,T;L^2(\Omega))$ and $h \in H^2(0,T;H^2(\Omega))$, with the compatibility conditions $h(x,0) = y^0(x)$ and $\partial_t h(x,0) = y^1(x)$ for all $x\in\partial\Omega$ then 
 $\partial_t y[p]$ solution of \eqref{WE} belongs to the space $C^0([0,T]; H^1(\Omega))\cap C^1([0,T]; L^2(\Omega))$, with estimates depending only on $m$ and the norms of $(y^0, y^1),\, g, \, h$ in these spaces. Therefore, due to Sobolev's imbedding, $\partial_t y[p] \in L^2 (0,T;L^{\infty}(\Omega))$. 
 \end{rem}

\begin{rem}\label{RemReg}
	The condition on $y^0$ requires compatibility conditions for $y^0$ and $h$ on $\partial\Omega\times \{0\}$. In particular,
$|h(x,0)|\geq r> 0,~x\in\partial\Omega$ must be satisfied since $|y^0(x)|\geq r> 0,~ae$ in $\overline{\Omega}$.
 \end{rem}

The method of proof of Theorem \ref{TWE} is based on a global Carleman estimate and is very close to the approach of \cite{ImYamCom01}, that concerns the wave equation with Neumann boundary condition and Dirichlet observation for the inverse problem of retrieving a potential. Actually, it also closely follows the approach of \cite{Yam99} but this work requires less regularity conditions on $y$. 

The use of Carleman estimate to prove uniqueness in inverse problems was introduced in \cite{BuKli81} by A. L. Bukhgeim and M. V. Klibanov.  Concerning inverse problems for hyperbolic equations with a single observation, we can refer to  \cite{PuelYam96}, \cite{PuelYam97}, where the method relies on uniqueness results obtained by local Carleman estimates  (see e.g. \cite{Im02}, \cite{Isakov}) and compactness-uniqueness arguments based on observability inequalities (see also \cite{Zhang00}). Related references \cite{ImYamCom01} and \cite{ImYamIP01} also use  Carleman estimates, but rather consider the case of interior or Dirichlet boundary data observation. 


\section{Classical Global Carleman Estimate}
In this section, we are interested in proving a global Carleman estimate for the wave operator. One can find (local) Carleman estimates for regular functions with compact support in \cite{Carleman}, \cite{HormanderLPDE} and in \cite{Yam99}.\\

Let us define the usual wave operator $L$ by 
$$
	Lv = \partial_t^2 v - \Delta v
$$
We consider a function $v\in L^2(-T,T;H^1_0(\Omega))$ such that $Lv \in L^2(-T,T;L^2(\Omega))$, 
and  satisfying $v(\pm T) = 0$, $\partial_t v(\pm T) = 0$ on $\Omega$.

Let us now define, for $x_{0} \not\in\Omega$, $\lambda>0$ and $\beta \in (0,1)$, and 
for all $(x,t) \in \Omega \times (-T,T)$:
\begin{equation}\label{weights}
	\psi(x,t) = |x-x_0|^2-\beta t^2+C_0 \quad \text{ and } \quad
	\varphi(x,t) = e^{\lambda \psi(x,t)}
\end{equation}

where $C_0>0$ is chosen so that $\psi\geq 1$ on $\Omega\times(-T,T)$. We also set, for $s>0$,
$$
	w(x,t)=v(x,t)e^{s\varphi(x,t)}.
$$
We define the operator $L_p$ for $p\in L^\infty_{\leq m} (\Omega)$ by
$$
L_pv = \partial_t^2 v - \Delta v + pv
$$
that satisfies $L_pv \in L^2(-T,T;L^2(\Omega))$ if $Lv \in L^2(-T,T;L^2(\Omega))$.\\

Let us first introduce the Carleman estimate we will prove by the formal calculation of 
$$
	Pw = e^{s\varphi}L(e^{-s\varphi} w).
$$
We have easily
\begin{eqnarray*}
Pw &=& \partial_t^2 w-2s\lambda\varphi (\partial_t w\partial_t \psi -\nabla w.\nabla\psi)+s^2\lambda^2\varphi^2w(|\partial_t \psi|^2-|\nabla\psi|^2)
-\Delta w\\
&&- s \lambda\varphi w(\partial_t^2 \psi - \Delta \psi )-s\lambda^2 \varphi w(|\partial_t \psi|^2-|\nabla\psi|^2)
\end{eqnarray*}
and if we set 
\begin{equation}\label{P}
	\begin{aligned}
	P_{1}w&=\partial_t^2 w-\Delta w+s^2\lambda^2\varphi^2w(|\partial_t \psi|^2-|\nabla\psi|^2)\\
	P_{2}w&=(\alpha - 1)s \lambda\varphi w(\partial_t^2 \psi - \Delta \psi )
		-s\lambda^2 \varphi w(|\partial_t \psi|^2-|\nabla\psi|^2)
		-2s\lambda\varphi (\partial_t w\partial_t \psi -\nabla w.\nabla\psi)\\
	Rw&=-\alpha s \lambda\varphi w(\partial_t^2 \psi - \Delta \psi )
	\end{aligned}
\end{equation}
with $\alpha$ chosen later  such that 
$ \frac {2\beta}{\beta+n} < \alpha < \frac {2}{\beta+n}$ (see \eqref{encadrement}),
we get $$P_{1}w+P_{2}w=Pw-Rw.$$

Let us now give a global (meaning ``up to the boundary'') Carleman inequality, following Imanuvilov's method \cite{Im02}.
One can read other versions of global Carleman estimates for hyperbolic equations in \cite{Zhang00} and \cite{Tataru96}.

\begin{theo}\label{CarlemanLucie}
Let us suppose that there exists $x_{0} \not\in\Omega$ such that 
$$
\Gamma_{0}\supset \{x\in\partial\Omega ; (x-x_{0})\cdot\nu(x)\geq 0\}
$$
Then for every $m>0$, there exists $\lambda_{0}>0$, $s_{0}>0$
and a constant $M=M(s_{0},\lambda_{0},T,m,\Omega,\beta,x_{0})$ 
such that for all $p\in L^\infty_{\leq m} (\Omega)$, and for all $\lambda>\lambda_{0}$, $s>s_{0}$:
\begin{gather}
	s\lambda \int_{-T}^{T} \int_{\Omega} e^{2s\varphi}(|\partial_t v|^2 + |\nabla v |^2)\,dxdt 
	+ s^3\lambda^3\int_{-T}^{T} \int_{\Omega} e^{2s\varphi}|v|^2\,dxdt\nonumber\\
	+ \int_{-T}^{T} \int_{\Omega} |P_1 (e^{s\varphi}v)|^2 \,dxdt + \int_{-T}^{T} \int_{\Omega} |P_2( e^{s\varphi}v)|^2\,dxdt\label{CarlemClas}\\
	\leq M\int_{-T}^{T} \int_{\Omega} e^{2s\varphi}| L_pv|^2\,dxdt 
	+ Ms\lambda \int_{-T}^{T} \int_{\Gamma_{0}} \varphi e^{2s\varphi} \left|\partial_\nu v\right|^2\,d\sigma dt.\nonumber
\end{gather}
for all $v\in H^1(-T,T;H^1_0(\Omega))$ satisfying
$$
	\left\lbrace\begin{array}{l}
		Lv \in L^2(\Omega\times (-T,T)),\\
		v(x,\pm T) = \partial_tv(x,\pm T)=0,~\forall x\in\Omega
	\end{array}\right.
$$
\end{theo}
We do not give here extensive references about Carleman estimates for hyperbolic equations, but one can have a look at \cite{Im02} and find references therein.
\begin{rem}
Estimate \eqref{CarlemClas} is uniform in $p$ when $p$ is in a bounded subset of $L^\infty_{\leq m}(\Omega)$.
Moreover, this Carleman estimate is proved for any arbitrary time $T$.  
\end{rem}

\noindent {\bf Proof.}

We will first prove estimate \eqref{CarlemClas} with $Lv$ in the right hand side instead of $L_p v$. One will see at the end that the result hold as well for $L_p$ since $p\in L^\infty(\Omega\times(-T,T))$.\\

As we began to write, for $w = e^{s\varphi} v$ we have $Pw = e^{s\varphi}Lv$ and 
\begin{gather}
	\int_{-T}^{T} \int_{\Omega}\left(| P_{1}w| ^2 + | P_{2}w| ^2\right)\,dxdt
	+~2\int_{-T}^{T} \int_{\Omega}P_{1}wP_{2}w\,dxdt \nonumber\\
	= \int_{-T}^{T} \int_{\Omega}| Pw-Rw|^2\,dxdt\label{PP}
\end{gather}
\bigskip
We will calculate and obtain a lower bound for 
$$
	\int_{-T}^{T} \int_{\Omega}P_{1}wP_{2}w\,dxdt.
$$ 

The main goal of the proof will be indeed to minimize this cross-term by positive and dominant terms looking similar to the one of the left hand side of \eqref{CarlemClas} and negative boundary terms that will be moved to the right hand side of the estimate. In the sake of clarity, we will devide the proof in several steps.\\

\noindent {\textbf{Step 1. Explicit calculations} 

We set $\ds\left\langle P_1w,P_2w\right\rangle_{L^2(\Omega\times(-T,T))}=\sum_{i,k=0}^n I_{i,k}$ where $I_{i,k}$ is the integral of the product of the $i$th-term in $P_1w$ and the $k$th-term in $P_2w$. 
We mainly use integrations by parts and the properties of $w$ such as $w(\pm T) = 0$, $\partial_tw(\pm T) = 0$ in $\Omega$ and $w=0$ on $\partial\Omega\times (-T,T)$.\\

 Integrations by part in time give easily, since $\partial_t\Delta \psi = 0$,
\begin{eqnarray*}
	I_{11}&=&\int_{-T}^{T} \int_{\Omega}\partial_t^2 w((\alpha - 1)s \lambda\varphi w(\partial_t^2 \psi - \Delta \psi ))\,dxdt\\
	&=&(1-\alpha )s\lambda \int_{-T}^{T} \int_{\Omega} \varphi  | \partial_t w | ^2 (\partial_t^2 \psi - \Delta \psi ) \,dxdt\\
	&&- \dfrac {(1-\alpha )}2 s\lambda^2 \int_{-T}^{T} \int_{\Omega} \varphi | w | ^2 \partial_t^2 \psi 	(\partial_t^2 \psi - \Delta \psi ) \,dxdt\\
	&&- \dfrac {(1-\alpha )}2 s\lambda^3 \int_{-T}^{T} \int_{\Omega} \varphi | w | ^2  |\partial_t \psi|^2(\partial_t^2 \psi - \Delta \psi ) \,dxdt.
\end{eqnarray*}
In the same manner, since $\partial_t \nabla \psi = \bf 0$, one has
\begin{eqnarray*}
	I_{12}&=&\int_{-T}^{T} \int_{\Omega}\partial_t^2 w(-s\lambda^2 \varphi w(|\partial_t \psi|^2-|\nabla\psi|^2))\,dxdt\\
	&=&s \lambda ^2 \int_{-T}^{T} \int_{\Omega} \varphi |\partial_t w|^2(|\partial_t \psi|^2-|\nabla\psi|^2)\,dxdt
	-s\lambda^2 \int_{-T}^{T} \int_{\Omega} \varphi  |w|^2 |\partial_t^2 \psi|^2\,dxdt\\
	&&-(2+\dfrac 12)s\lambda^3 \int_{-T}^{T} \int_{\Omega} \varphi |w|^2 |\partial_t \psi|^2\partial_t^2 \psi\,dxdt
	+\dfrac {s\lambda^3}2 \int_{-T}^{T} \int_{\Omega} \varphi |w|^2 |\nabla\psi|^2\partial_t^2 \psi\,dxdt\\
	&&-\dfrac {s\lambda^4}2 \int_{-T}^{T} \int_{\Omega} \varphi  |w|^2 |\partial_t \psi|^2 (|\partial_t \psi|^2-|\nabla\psi|^2)\,dxdt
\end{eqnarray*}
and using also integrations by part in space variable, we get
\begin{eqnarray*}
	I_{13}&=&\int_{-T}^{T} \int_{\Omega}\partial_t^2 w(-2s\lambda\varphi (\partial_t w\partial_t \psi -\nabla w.\nabla\psi))\,dxdt\\
	&=& s\lambda \int_{-T}^{T} \int_{\Omega} \varphi | \partial_t w | ^2 \partial_t^2 \psi \,dxdt
	+s\lambda^2 \int_{-T}^{T} \int_{\Omega} \varphi | \partial_t w | ^2 |\partial_t \psi|^2 \,dxdt\\
	&&+ s\lambda \int_{-T}^{T} \int_{\Omega} \varphi | \partial_t w | ^2 \Delta \psi  \,dxdt
	+ s\lambda^2 \int_{-T}^{T} \int_{\Omega} \varphi | \partial_t w | ^2  |\nabla\psi|^2 \,dxdt\\
	&&-2s\lambda^2 \int_{-T}^{T} \int_{\Omega} \varphi \partial_t w\, \partial_t \psi \nabla w.\nabla\psi \,dxdt.
\end{eqnarray*}

We compute in the same way
\begin{eqnarray*}
	I_{21}&=&\int_{-T}^{T} \int_{\Omega}-\Delta w ((\alpha - 1)s \lambda\varphi w(\partial_t^2 \psi - \Delta \psi ))\,dxdt\\
	&=& -(1-\alpha ) s \lambda \int_{-T}^{T} \int_{\Omega} \varphi |\nabla w|^2 (\partial_t^2 \psi - \Delta \psi ) \,dxdt\\
	&&+ \dfrac {(1-\alpha )}2 s \lambda^2 \int_{-T}^{T} \int_{\Omega} \varphi |w|^2 \Delta \psi (\partial_t^2 \psi - \Delta \psi )\,dxdt\\
	&&+ \dfrac {(1-\alpha )}2 s \lambda^3 \int_{-T}^{T} \int_{\Omega} \varphi |w|^2  |\nabla\psi|^2 (\partial_t^2 \psi - \Delta \psi )\,dxdt
\end{eqnarray*}
and 
\begin{eqnarray*}
	I_{22}&=&\int_{-T}^{T} \int_{\Omega}-\Delta w(-s\lambda^2\varphi w(|\partial_t \psi|^2-|\nabla\psi|^2))\,dxdt\\
	&=&-s \lambda^2 \int_{-T}^{T} \int_{\Omega} \varphi |\nabla w|^2(|\partial_t \psi|^2-|\nabla\psi|^2) \,dxdt\\
	&&-\dfrac{s \lambda^2}2 \int_{-T}^{T} \int_{\Omega} \varphi |w|^2 \Delta(|\nabla \psi|^2) \,dxdt\\
	&&+\dfrac {s \lambda^3}2 \int_{-T}^{T} \int_{\Omega} \varphi  |w|^2\Delta \psi (|\partial_t \psi|^2-|\nabla\psi|^2)\,dxdt\\
	&&+\dfrac {s \lambda^4}2 \int_{-T}^{T} \int_{\Omega} \varphi  |w|^2 |\nabla\psi|^2(|\partial_t \psi|^2-|\nabla\psi|^2)\,dxdt\\
	&&-s \lambda^3\int_{-T}^{T} \int_{\Omega} \varphi |w|^2 \nabla\psi \cdot \nabla(|\nabla\psi|^2) \,dxdt
\end{eqnarray*}
Using the fact that $w|_{\partial\Omega\times (-T,T)}=0$, we have, on $\partial\Omega\times (-T,T)$, $\nabla w = (\partial_\nu w) \nu$ that give $|w|^2 = |\partial_\nu w|^2$. Therefore, we obtain
\begin{eqnarray*}
	I_{23}&=&\int_{-T}^{T} \int_{\Omega}-\Delta w(-2s\lambda\varphi (\partial_t w\partial_t \psi -\nabla w \cdot \nabla\psi))\,dxdt\\
	&=& s\lambda \int_{-T}^{T} \int_{\Omega} \varphi |\nabla w|^2 (\partial_t^2 \psi - \Delta \psi )\,dxdt
	+2s\lambda^2 \int_{-T}^{T} \int_{\Omega} \varphi | \nabla \psi\cdot\nabla w |^2 \,dxdt\\
	&&-2s\lambda^2 \int_{-T}^{T} \int_{\Omega} \varphi\partial_t w\, \partial_t \psi \nabla w\cdot\nabla\psi \,dxdt
	+s\lambda^2 \int_{-T}^{T} \int_{\Omega} \varphi | \nabla w |^2  (|\partial_t \psi|^2-|\nabla\psi|^2) \,dxdt\\
	&&- s\lambda \int_{-T}^{T} \int_{\partial\Omega} \varphi | \partial_\nu w|^2 \nabla \psi \cdot \nu \,d\sigma dt
	+2s\lambda \int_{-T}^{T} \int_{\Omega} \varphi  D^2 \psi |\nabla w|^2\,dxdt
\end{eqnarray*}
since $D^2 \psi$ is symmetric.

One easily write
\begin{eqnarray*}
	I_{31}&=&\int_{-T}^{T} \int_{\Omega}  s^2\lambda^2\varphi^2w(|\partial_t \psi|^2-|\nabla\psi|^2)
	((\alpha - 1)s \lambda\varphi w(\partial_t^2 \psi - \Delta \psi ))\,dxdt\\
	&=&(\alpha -1) s^3\lambda^3 \int_{-T}^{T} \int_{\Omega} \varphi^3 |w|^2(\partial_t^2 \psi - \Delta \psi ) (|\partial_t \psi|^2-|\nabla\psi|^2)\,dxdt
\end{eqnarray*}
and
\begin{eqnarray*}
	I_{32}&=&\int_{-T}^{T} \int_{\Omega}   s^2\lambda^2\varphi^2w(|\partial_t \psi|^2-|\nabla\psi|^2)
	(-s\lambda^2 \varphi w(|\partial_t \psi|^2-|\nabla\psi|^2))\,dxdt\\
	&=&-s^3\lambda^4 \int_{-T}^{T} \int_{\Omega}\varphi^3  |w|^2 (|\partial_t \psi|^2-|\nabla\psi|^2)^2\,dxdt.
\end{eqnarray*}
Finally, some integrations by part enable to obtain, since $\nabla(|\partial_t \psi|^2) = \bf 0$,
\begin{eqnarray*}
	I_{33}&=&\int_{-T}^{T} \int_{\Omega}  s^2\lambda^2\varphi^2w(|\partial_t \psi|^2-|\nabla\psi|^2) 
	(-2s\lambda  \varphi (\partial_t w \partial_t \psi -\nabla w.\nabla\psi)) \,dxdt\\
	&=&s^3\lambda ^3 \int_{-T}^{T} \int_{\Omega}\varphi^3 |w|^2 (\partial_t^2 \psi - \Delta \psi ) (|\partial_t \psi|^2-|\nabla\psi|^2) \,dxdt\\
	&&+s^3\lambda ^3 \int_{-T}^{T} \int_{\Omega}\varphi^3  |w^2| 
		(2\partial_t^2 \psi|\partial_t \psi|^2 + \nabla\psi \cdot \nabla(|\nabla\psi|^2))\,dxdt\\
	&&+3s^3\lambda^4 \int_{-T}^{T} \int_{\Omega}\varphi^3  |w|^2  (|\partial_t \psi|^2-|\nabla\psi|^2)^2\,dxdt.
\end{eqnarray*}

Gathering all the terms that have been computed, we get 
\begin{equation}\label{P1P2}
	\begin{aligned}
	\int_{-T}^{T} \int_{\Omega} P_{1}&wP_{2}w\,dxdt\\
	=& ~2s\lambda \int_{-T}^{T} \int_{\Omega} \varphi | \partial_t w | ^2 \partial_t^2 \psi \,dxdt
	- \alpha s\lambda \int_{-T}^{T} \int_{\Omega} \varphi | \partial_t w | ^2 (\partial_t^2 \psi - \Delta \psi ) \,dxdt\\
	&+~2s\lambda^2 \int_{-T}^{T} \int_{\Omega}\varphi \left( | \partial_t w | ^2 |\partial_t \psi|^2
	-2\partial_t w\, \partial_t \psi  \nabla w\cdot \nabla\psi + | \nabla \phi\cdot\nabla w |^2 \right) \,dxdt\\
	&+~2s\lambda \int_{-T}^{T} \int_{\Omega} \varphi  D^2 \psi |\nabla w|^2\,dxdt
	+\alpha s \lambda \int_{-T}^{T} \int_{\Omega} \varphi |\nabla w|^2 (\partial_t^2 \psi - \Delta \psi ) \,dxdt\\
	&-~s\lambda \int_{-T}^{T} \int_{\partial\Omega}\varphi \left| \partial_\nu w\right |^2 \nabla\psi\cdot\nu(x)\,d\sigma dt\\
	&+~2s^3\lambda^4 \int_{-T}^{T} \int_{\Omega} \varphi^3 |w|^2 (|\partial_t \psi|^2-|\nabla\psi|^2)^2\,dxdt\\
	&+~s^3\lambda ^3 \int_{-T}^{T} \int_{\Omega}\varphi^3  |w^2| 
		(2\partial_t^2 \psi|\partial_t \psi|^2 + \nabla\psi \cdot \nabla(|\nabla\psi|^2) )\,dxdt\\
	&+~\alpha s^3\lambda^3 \int_{-T}^{T} \int_{\Omega}\varphi^3  |w|^2
		(\partial_t^2 \psi - \Delta \psi ) (|\partial_t \psi|^2-|\nabla\psi|^2)\,dxdt\\
	&+~ X_1
	\end{aligned}
\end{equation}
where $X_1$ satisfies, using the regularity of $\psi$ and the fact that  $\psi \geq 1$ implies $\lambda \leq e^{\lambda \psi} = \varphi$,
\begin{equation}\label{X1}
|X_1| \leq M s\lambda^3 \int_{-T}^{T} \int_{\Omega} \varphi^3 |w|^2  \,dxdt .
\end{equation}
Here and in the sequel, $M>0$ corresponds to a generic constant depending at least on $\Omega$ and $T$ but independant of $s$ and $\lambda$.\\

\noindent {\textbf{Step 2. Bounding each term from below} 

On the one hand, one can notice that 
\begin{multline}\label{00}
	2s\lambda^2 \int_{-T}^{T} \int_{\Omega}\varphi \left( | \partial_t w | ^2 |\partial_t \psi|^2
		-2\partial_t w \partial_t \psi  \nabla w \cdot \nabla\psi + | \nabla \psi\cdot\nabla w |^2 \right)\,dxdt\\
	=2s\lambda^2 \int_{-T}^{T} \int_{\Omega}\varphi \left( \partial_t w \partial_t \psi- \nabla w \cdot \nabla\psi\right)^2 \,dxdt \geq 0.
\end{multline}
Besides, considering the terms in $s\lambda$ that must now give the dominant terms in $| \partial_t w|^2 $ and $|\nabla w|^2$
and thus have to be strictly positive, one can guess that we need
$$
	2 \partial_t^2 \psi -\alpha (\partial_t^2 \psi - \Delta\psi) > 0 \quad \tn{ and } \quad
	2D^2\psi +\alpha  (\partial_t^2 \psi - \Delta \psi ) > 0.
$$
This will constrain the value of the unspecified constant $\alpha>0$ in the definition \eqref{P} of $P_2 w$.
We get, by explicit computations, that $\beta \in (0,1)$, and $\alpha$ must satisfy
\begin{equation}\label{encadrement}
	 \dfrac{2\beta}{\beta +n} < \alpha <  \dfrac{2}{\beta +n}.
\end{equation}
As a direct consequence, we can write 
\begin{multline}
	\label{ordre1}
	2 s\lambda \int_{-T}^{T} \int_{\Omega} \varphi | \partial_t w|^2\partial_t^2 \psi \,dxdt 
	-\alpha s\lambda \int_{-T}^{T} \int_{\Omega} \varphi | \partial_t w|^2(\partial_t^2 \psi - \Delta\psi) \,dxdt\\
	+2s\lambda \int_{-T}^{T} \int_{\Omega} \varphi D^2\psi |\nabla w|^2 \,dxdt
	+\alpha s \lambda \int_{-T}^{T} \int_{\Omega} \varphi |\nabla w|^2 (\partial_t^2 \psi - \Delta \psi ) \,dxdt\\
	\geq M s\lambda \int_{-T}^{T} \int_{\Omega} \varphi | \partial_t w|^2\,dxdt 
	+ M s\lambda \int_{-T}^{T} \int_{\Omega} \varphi | \nabla w|^2 \,dxdt.
\end{multline}

On the other hand, we can observe that:
\begin{eqnarray*}
	&&2s^3\lambda^4 \int_{-T}^{T} \int_{\Omega}\varphi^3 |w|^2  (|\partial_t \psi|^2-|\nabla\psi|^2)^2\,dxdt\\
	&&+~s^3\lambda ^3 \int_{-T}^{T} \int_{\Omega}\varphi^3 |w^2| 
		(2\partial_t^2 \psi|\partial_t \psi|^2 + \nabla\psi \cdot \nabla(|\nabla\psi|^2) )\,dxdt\\
	&&+~\alpha s^3\lambda^3 \int_{-T}^{T} \int_{\Omega}\varphi^3 |w|^2(\partial_t^2 \psi - \Delta \psi )
		 (|\partial_t \psi|^2-|\nabla\psi|^2)\,dxdt\\
	&=& s^3\lambda ^3 \int_{-T}^{T} \int_{\Omega}\varphi^3 |w^2|  F_{\lambda}(\phi) \,dxdt
\end{eqnarray*}
where 
\begin{align*}
	F_{\lambda}(\phi)&= 
	2\lambda(|\partial_t \psi|^2-|\nabla\psi|^2)^2 +  (2\partial_t^2 \psi|\partial_t \psi|^2 + \nabla\psi \cdot \nabla(|\nabla\psi|^2) )
		+ \alpha (\partial_t^2 \psi - \Delta \psi ) (|\partial_t \psi|^2-|\nabla\psi|^2)\\
	&= 32 \lambda (\beta^2t^2 - |x-x_0|^2)^2 - 16(\beta^3t^2 - |x-x_0|^2) -8 \alpha(\beta + n)(\beta^2t^2 - |x-x_0|^2)\\
	&= 32 \lambda (\beta^2t^2 - |x-x_0|^2)^2 - 8(\alpha(\beta + n) +2\beta) (\beta^2t^2 - |x-x_0|^2) +16(1-\beta)|x-x_0|^2.
\end{align*}
Since $x_0\not\in\Omega$, we have $16(1-\beta)|x-x_0|^2\geq c_*>0$. 
Therefore, we are considering a polynome $P(X) \geq 32\lambda X^2 - 8\left( \alpha (\beta +n) + 2 \beta \right) X+c_*$ and taking $\lambda>0$ large enough, the minimum of $P$ will be strictly positive.
Consequently, 
\begin{equation}
	\label{ordre0}
	s^3\lambda ^3 \int_{-T}^{T} \int_{\Omega} \varphi^3|w^2| F_{\lambda}(\phi) \,dxdt\\
	\geq M s^3\lambda^3 \int_{-T}^{T} \int_{\Omega} \varphi^3 | w|^2\,dxdt.
\end{equation}\\

Thus, plugging \eqref{00}, \eqref{ordre1} and \eqref{ordre0} in \eqref{P1P2}, we  obtain
\begin{multline*}
	\int_{-T}^{T} \int_{\Omega}P_{1}wP_{2}w\,dxdt  
	+ 2 s \lambda\int_{-T}^{T} 	\int_{\partial\Omega}\varphi \left| \partial_\nu w\right|^2 (x-x_0)\cdot\nu(x)\,d\sigma dt - X_1\\
	\geq M s\lambda \int_{-T}^{T} \int_{\Omega} \varphi\left(| \partial_t w| ^2+|\nabla w|^2\right)\,dxdt
		+ M s^3\lambda ^3\int_{-T}^{T} \int_{\Omega}\varphi^3 | w| ^2\,dxdt.
\end{multline*}
Since we also have
\begin{multline*}
	\int_{-T}^{T} \int_{\Omega}| Pw-Rw|^2\,dxdt 
		\leq 2\int_{-T}^{T} \int_{\Omega}| Pw|^2\,dxdt 
		+2\int_{-T}^{T} \int_{\Omega}| Rw|^2\,dxdt\\
	\leq M\int_{-T}^{T} \int_{\Omega}| Pw|^2\,dxdt 
	+Ms^2 \lambda^2\int_{-T}^{T} \int_{\Omega} \varphi^2 |w|^2 \,dxdt,
\end{multline*}
using \eqref{PP} and \eqref{X1}, we get 
\begin{gather*}
	s\lambda \int_{-T}^{T} \int_{\Omega}\left(| \partial_t w| ^2+|\nabla w|^2\right) \varphi\,dxdt
		+ s^3\lambda ^3\int_{-T}^{T} \int_{\Omega} | w| ^2\varphi^3\,dxdt\\
	+\int_{-T}^{T} \int_{\Omega}\left(| P_{1}w| ^2 + | P_{2}w| ^2\right)\,dxdt\\
	\leq M\int_{-T}^{T} \int_{\Omega}| Pw|^2\,dxdt 
		+M s \lambda\int_{-T}^{T} \int_{\partial\Omega}\varphi \left| \partial_\nu w\right|^2 
			(x-x_0)\cdot\nu(x)\,d\sigma dt \\
	+~M s\lambda ^3\int_{-T}^{T} \int_{\Omega}\varphi^3 | w| ^2\,dxdt  + Ms^2 \lambda^2\int_{-T}^{T} \int_{\Omega} \varphi^2 |w|^2 \,dxdt.
\end{gather*}
We take now $s_0$ large enough so that the terms of the last line (coming from  $X_1$ and $| Rw|^2$) are absorbed by the dominant term in $s^3\lambda^3  |w|^2\varphi^3$ as soon as $s>s_0$. Using also the condition \eqref{GCC} on $\Gamma_0$, we finally obtain for some positive constant $M=M(s_{0},\lambda_{0},m,\Omega,\beta,x_{0})$,
\begin{gather}
	s\lambda\int_{-T}^{T} \int_{\Omega} \varphi(| \partial_t w|^2+|\nabla w|^2)\,dxdt
	+ s^3\lambda ^3\int_{-T}^{T} \int_{\Omega} \varphi^3 |w|^2\,dxdt\nonumber\\
	+\int_{-T}^{T} \int_{\Omega}|P_1 w|^2\,dxdt+\int_{-T}^{T} \int_{\Omega}|P_2 w|^2\,dxdt\label{CarlemW}\\
	\leq M\int_{-T}^{T} \int_{\Omega}| Pw|^2\,dxdt 
	+Ms\lambda \int_{-T}^{T} \int_{\Gamma_0}  \varphi\left|\partial_\nu w\right|^2\,d\sigma dt \nonumber
\end{gather}
$\forall s>s_0$, $\forall \lambda > \lambda_0$.\\

\noindent {\textbf{Step 3. Return to the variable $v$} 

Using that $w=ve^{s\varphi}$ 
gives for all $x\in \Omega$ and $t\in(-T,T)$
\begin{eqnarray*}
	&&e^{2s\varphi}|\partial_t v|^2 \leq 2|\partial_t w|^2+2s^2|\partial_t\varphi|^2|w|^2,\\
	&&e^{2s\varphi}|\nabla v|^2 \leq 2|\nabla w|^2+2s^2|\nabla\varphi|^2|w|^2,\\
	&&e^{2s\varphi}\left|\partial_\nu v\right|^2= \left|\partial_\nu w\right|^2 \tn{ on }\partial\Omega
\end{eqnarray*}
and that by construction $P w = e^{s\varphi} Lv$, 
we can go back to the variable $v$ in \eqref{CarlemW} and obtain that there exists some positive constant $M=M(s_{0},\lambda_{0},m,\Omega,\beta,x_{0})$ such that for all $ s>s_0$ and $ \lambda > \lambda_0$,
\begin{gather}
	s\lambda \int_{-T}^{T} \int_{\Omega} e^{2s\varphi}(|\partial_t v|^2 + |\nabla v |^2)\,dxdt 
	+ s^3\lambda^3\int_{-T}^{T} \int_{\Omega} e^{2s\varphi}|v|^2\,dxdt\nonumber \\
	+ \int_{-T}^{T} \int_{\Omega} |P_1 (e^{s\varphi} v)|^2 \,dxdt + \int_{-T}^{T} \int_{\Omega} |P_2 (e^{s\varphi} v)|^2\,dxdt \label{CarlemL}\\
	\leq M\int_{-T}^{T} \int_{\Omega} e^{2s\varphi}| Lv|^2\,dxdt 
	+ Ms\lambda \int_{-T}^{T} \int_{\Gamma_0} \varphi e^{2s\varphi} \left|\partial_\nu v \right|^2 \,d\sigma dt. \nonumber
\end{gather}
It concludes the proof of a Carleman estimate for the operator $L=\partial_t^2 - \Delta$.\\

\noindent {\textbf{Step 4. Wave operator with potential} 

The Carleman estimate \eqref{CarlemClas} for the operator $L_p=\partial_t^2 - \Delta + p$ with $p\in L^\infty_{\leq m} (\Omega)$ is a direct consequence of \eqref{CarlemL} noticing that on $\Omega\times(-T,T)$,
$$
	|Lv|^2 \leq 2 |L_pv|^2 + 2 \|p\|_{L^\infty_{\leq m} (\Omega)}|v|^2 
	\leq 2 |L_pv|^2 + 2 m |v|^2.
$$
Indeed, choosing $s_0$ (or $\lambda_0$) large enough, one can absorb the term
$$
	2M m\int_{-T}^{T} \int_{\Omega} e^{2s\varphi}|v|^2\,dxdt
$$
in the left hand side of \eqref{CarlemL} and obtain \eqref{CarlemClas} with slightly different constants.
This ends the proof of Theorem~\ref{CarlemanLucie}.

\section{Stability of the inverse problem}

Before giving the proof of Theorem~\ref{TWE}, we will begin this section by a stability theorem for the linear inverse source problem stated after equation~\eqref{WE2}. The following answer is obtained using the global Carleman estimate given above.
\begin{theo}\label{Thm2}
Let $m>0$, $K>0$ and $r>0$. Let $q$ belong to  $L^\infty_{\leq m}(\Omega)$. 
Assume that $f\in L^2(\Omega)$ and $R\in H^1(0,T;L^{\infty}(\Omega))$ with
$$|| R||_{H^1 (0,T;L^{\infty}(\Omega))} \leq K $$
and 
\begin{equation} \label{RR}
\inf_{x\in\Omega} |R(x,0)| \geq r.
\end{equation}
Under the hypothesis that $\Gamma_0$  and $T$ satisfy the geometrical conditions 
$$
\exists x_{0}\not\in\Omega\quad \textit{such that} \quad \Gamma_{0}\supset \{x\in\partial\Omega ; (x-x_{0})\cdot \nu(x)\geq 0\}
\quad \textit{and }\quad T> \sup_{x\in \Omega} |x-x_0|
$$
if $u[f]$ is the solution of equation~\eqref{WE2}: 
$$
	\left\lbrace\begin{array}{ll}
		\partial_t^2u(x,t)-\Delta u(x,t)+q(x)u(x,t)=f(x)R(x,t),&\quad x\in\Omega,~t\in (0,T)\\
		u(x,t)=0,&\quad x \in\partial\Omega,~t\in (0,T)\\
		u(x,0)=0,~\partial_t u(x,0)= 0,&\quad x\in\Omega,
	\end{array} \right.
$$
then there exists a constant $C=C(T,\Omega,\Gamma_0,K,r)>0$ such that for all $f \in L^2(\Omega)$ :
\begin{equation}\label{stablin}
C^{-1}|| f||_{L^2(\Omega)}\leq   \left\|\dfrac{ \partial u[f]}{\partial\nu} \right\|_{H^1(0,T;L^2(\Gamma_0))} \leq C || f||_{L^2(\Omega)} .
\end{equation}
\end{theo}

\noindent {\bf Proof.}
We will apply the Carleman estimate given by \eqref{CarlemL} to $w=\chi \partial_t y$ where $\chi$ is a cutoff function to be detailed later. We divide the proof in several steps.\\

\noindent {\textbf{Step 1.} }
Let us first work on the equation satisfied by $z=\partial_t u$:
\begin{equation}\label{WE3}
	\left\lbrace\begin{array}{ll}
		\partial_t^2z(x,t)-\Delta z(x,t)+q(x)z(x,t)=f(x)\partial_t R(x,t),&\quad x\in\Omega,~t\in (0,T)\\
		z(x,t)=0,&\quad x \in\partial\Omega,~t\in (0,T)\\
		z(x,0)=0,\quad \partial_t z(x,0)= f(x) R(x,0),&\quad x\in\Omega.
	\end{array} \right.
\end{equation}
We gather in the following lemma the classical energy and trace estimates we will need in the sequel. 
\begin{lemme}\label{Energy}
Assume that $p\in L_{\leq m}^{\infty}(\Omega)$, $g \in L^1(0,T;L^2(\Omega))$, $u^0\in H_0^1(\Omega)$ and $u^1\in L^2(\Omega)$
and consider the classical wave equation 
\begin{equation}\label{e}
	\left\lbrace\begin{array}{ll}
		\partial_t^2v(x,t)-\Delta v(x,t)+q(x)v(x,t)=g(x,t),&\quad x\in\Omega,~t\in (0,T)\\
		v(x,t)=0,&\quad x \in\partial\Omega,~t\in (0,T)\\
		v(x,0)=v_0(x),\quad \partial_t v(x,0)= v_1(x)),&\quad x\in\Omega.
	\end{array} \right.
\end{equation}
The Cauchy problem is well-posed and equation \eqref{e} admits a unique weak solution 
$$
	v\in C([0,T],H_0^1(\Omega)) \cap C^1([0,T],L^2(\Omega))
$$ 
and  there exists a constant 
$C=C(\Omega,T,m)>0$ such that for all $t\in(0,T)$, the energy 
$E_v(t) = ||\partial_t v(t)||^2_{L^2(\Omega)}+||\nabla v(t)||^2_{L^2(\Omega)}$ of the system satisfies
\begin{equation}\label{estimeenergy}
	E_v(t) \leq C\left(||v_0||^2_{H_0^1(\Omega)} + ||v_1||^2_{L^2(\Omega)} +
		||g||^2_{L^1(0,T,L^2(\Omega))} \right).
\end{equation}
Moreover, the normal derivative $\partial_\nu v$ belongs to $ L^2(0,T;L^2(\partial\Omega)) $ and satisfies
\begin{equation}\label{hiddenregularity}
	\left\|\dfrac{ \partial v}{\partial\nu} \right\|_{L^2(0,T;L^2(\partial\Omega))}\leq C \left(||v_0||^2_{H_0^1(\Omega)} 
		+ ||v_1||^2_{L^2(\Omega)} + ||g||^2_{L^1(0,T,L^2(\Omega))} \right).
\end{equation}
\end{lemme}

This result is very classical and we  refer to \cite{LionsMagenesBook} (Chapter 3) for the proof of the existence and uniqueness of solution to equation~\eqref{e}. Estimate \eqref{estimeenergy} can formally be deduced from multiplication of \eqref{e} by $v_t$, and the integration of this equality on $(0,T)\times \Omega$, using some integrations by parts. Concerning estimate~\eqref{hiddenregularity}, we refer to \cite{Lions} (Chapter 1). This estimate is a hidden regularity result which can be obtained by multipliers technique. \\

We can apply this lemma to equation \eqref{WE3} since $f\in L^2(\Omega)$ and $R\in H^1(0,T;L^{\infty}(\Omega))$. 
We denote the energy of the system by
$$
	E_z (t) = ||z(t)||^2_{H^1_0(\Omega)} + ||\partial_t z(t)||^2_{L^2(\Omega)}
$$ 
and we get 
$
	z \in C([0,T],H_0^1(\Omega)) \cap C^1([0,T],L^2(\Omega))
$ with, for all $t \in (0,T)$:
\begin{eqnarray}
	E_z (t) &\leq& C E_z(0) + C \| f\partial_t R \|^2_{L^1(0,T;L^2(\Omega))} \nonumber\\
	&\leq& C||f||^2_{L^2(\Omega)}\left( ||R(0)||^2_{L^{\infty}(\Omega)} + ||R||^2_{H^1(0,T;L^{\infty}(\Omega))}\right),\label{energy}
\end{eqnarray} 
and 
$
\partial_\nu z \in L^2(0,T;L^2(\partial\Omega))
$ with
$$
	\left\|\dfrac{ \partial z}{\partial\nu} \right\|_{L^2(0,T;L^2(\partial\Omega))}
	\leq C ||f||^2_{L^2(\Omega)}\left( ||R(0)||^2_{L^{\infty}(\Omega)} + ||R||^2_{H^1(0,T;L^{\infty}(\Omega))}\right).
$$
This last estimate gives at the end $\partial_\nu u \in H^1(0,T;L^2(\partial\Omega))$ 
and proves the right hand side of the two sided estimate~\eqref{stablin}. It thus gives a meaning to the measurement of the flux of the solution $u[f]$ we make in our inverse problem. \\

\noindent {\textbf{Step 2.}} Let us now extend the problem \eqref{WE3} on $(-T,T)$ in an even way, setting $z(x,t)=z(x,-t)$ for all $(x,t)\in\Omega\times(-T,0)$. We also extend $\partial_t R$ on an even way and keep the same notations for the new problem. Therefore, we have
$$z\in C([-T,T],H_0^1(\Omega)) \cap C^1([-T,T],L^2(\Omega))\tn{ and }\partial_t R \in L^2(-T,T;L^{\infty}(\Omega)).$$

The goal we have in mind is to apply the Carleman estimate to a solution of the wave equation in order to be able to estimate a weighted norm of $f$ by 
a weighted norm of the measurement $\partial_\nu u$. The solution $z$ could be a good candidate but is not equal to zero at time $\pm T$. This motivates the use of a cutoff function in time, that is calibrated to work well with the Carleman weight function $\psi$. 

We recall that for $\beta \in (0,1)$, $\psi$ and $\varphi$ are defined by \eqref{weights}.
Since $T$ satisfies \eqref{GCCT}, one can choose $\beta$ such that $\beta T^2 >  \sup_{x\in \Omega} |x-x_0|^2$ so that for all $x\in\Omega$, 
$\psi(x,\pm T)=|x-x_0|^2-\beta T^2 +C_0 < C_0$.
We then choose $\eta > 0$ such that
$$
	\psi(x,t) \leq C_0,~\quad \forall (x,t)\in\Omega\times[-T,-T+\eta] \cup [T-\eta,T].
$$
Besides, one observes that
$$
\psi(x,0)=|x-x_0|^2+C_0 > C_0,~ \quad\forall x\in\Omega.
$$ 
Therefore, for all $(x,t)\in\Omega\times[-T,-T+\eta] \cup [T-\eta,T]$,
\begin{equation}\label{psi}
	\varphi(x,t) \leq e^{\lambda C_0}< e^{\lambda \psi(x,0)}=\varphi(x,0)
\end{equation}
but we also have, for all $(x,t)\in\Omega\times[-T,T]$,
\begin{equation}\label{psi2}
	\varphi(x,t) \leq \varphi(x,0).
\end{equation}

Using the parameter $\eta$ introduced here, we define the cut-off function $\chi\in C^{\infty}(\mathbb{R};[0,1])$ such that
\begin{equation}\label{chi}
\left\lbrace
\begin{array}{ll}
\chi(\pm T) = \chi ' (\pm T) = 0&\\
\chi(t) =1&\forall t\in [-T+\eta,T-\eta]
\end{array} \right.
\end{equation}
and  we set $~v=\chi z~$ that satisfies the following equation:
\begin{equation}\label{WE4}
	\left\lbrace\begin{array}{ll}
		\partial_t^2v-\Delta v+qv= \chi f\partial_t R+ \chi''z+ 2 \chi' \partial_{t}z,&\quad \Omega\times(0,T)\\
		v(x,t)=0,&\quad x \in\partial\Omega,~t\in (0,T)\\
		v(x,0)=0,\quad \partial_t v(x,0)= f(x) R(x,0),&\quad x\in\Omega.
	\end{array} \right.
\end{equation}

Henceforth, $M>0$ will correspond to a generic constant depending on $s_0, \lambda_0,T,\Omega$, $\Gamma_0$, $\beta$, $\chi$, $r$, $K$, and $\eta$  but independent of $s>s_0$.\\

\noindent {\textbf{Step 3.} }
We use now the same notations as in the proof of the Carleman estimate. For $v$ solution of \eqref{WE4}, we set
$$w=e^{s\varphi} v~,~~P_{1}w=\partial_t^2 w-\Delta w+s^2\lambda^2\varphi^2w(|\partial_t \psi|^2-|\nabla\psi|^2),$$
and we then have for all $x\in \Omega$,
$w(x,\pm T)  =  \partial_t w(x,\pm T) = 0 $ and $ w(x,0)  =  0$ and 
for all $(x,t)\in \partial\Omega\times(-T,T)$, $w(x,t) = 0$.

Inspired from an idea of O. Yu. Imanuvilov and M. Yamamoto \cite{ImYamCom01}, we consider the integral 
$$\int_{-T}^0\int_{\Omega}P_1w(x,t)\,\partial_t w(x,t) \,dxdt.$$
On the one hand, using the properties of $w$, we can make the following calculation:
\begin{eqnarray*}
	&&\int_{-T}^0\int_{\Omega}P_1w\,\partial_t w \,dxdt\\
	& = &\int_{-T}^0\int_{\Omega} \left(\partial_t^2 w-\Delta w+s^2\lambda^2\varphi^2w(|\partial_t \psi|^2-|\nabla\psi|^2)\right)\partial_t w  \,dxdt\\
	& = & \dfrac 12 \int_{\Omega} |\partial_t w(0)|^2 \,dx
		- \dfrac {s^2\lambda^2}2  \int_{-T}^0\int_{\Omega} |w|^2 \partial_t \left(\varphi^2(|\partial_t \psi|^2-|\nabla\psi|^2)\right) \,dxdt\\
	& = & \dfrac 12 \int_{\Omega} |\partial_t w(0)|^2 \,dx
		 - 2 s^2\lambda^2  \int_{-T}^0\int_{\Omega}  |w|^2 \partial_t \left(\varphi^2(\beta^2 t^2-|x-x_0|^2)\right)\,dxdt\\
	& = & \dfrac 12\int_{\Omega} e^{2s\varphi(0)}|fR(0)|^2\,dx 
		+ 8 s^2\lambda^3  \int_{-T}^0\int_{\Omega}  |w|^2\varphi^2 \beta t(\beta^2 t^2-|x-x_0|^2) \,dxdt\\
	& &- ~4 s^2\lambda^2  \int_{-T}^0\int_{\Omega}  
|w|^2\varphi^2\beta^2 t.
\end{eqnarray*}
On the other hand, from this equality and a Cauchy-Schwarz estimate,
\begin{eqnarray}
	 \int_{\Omega} e^{2s\varphi(0)}|fR(0)|^2 \,dx
	& = &2 \int_{-T}^0\int_{\Omega}P_1w\, \partial_t w \,dxdt
		+ 8 s^2\lambda^2  \int_{-T}^0\int_{\Omega}  w^2\varphi^2\beta^2 t\,dxdt \nonumber\\
	&&- 16 s^2\lambda^3  \int_{-T}^0\int_{\Omega}  |w|^2\varphi^2 \beta t(\beta^2 t^2-|x-x_0|^2) \,dxdt \nonumber\\
	&\leq& 2 \left(\int_{-T}^0\int_\Omega\left| P_1w\right|^2\,dxdt \right)^{\frac 12}
		\left(\int_{-T}^0\int_ \Omega |\partial_t w|^2 \,dxdt\right)^{\frac 12} \nonumber\\
	&&+ M s^2\lambda^3 \int_{-T}^0\int_ \Omega  |w|^2\varphi^2 \left|\beta t(\beta^2 t^2-|x-x_0|^2)\right| \,dxdt\label{rhsw}.
\end{eqnarray}
Using now the `intermediate'' Carleman estimate~\eqref{CarlemW}  for a fixed $\lambda>\lambda_0$ not mentioned anymore, choosing $s$ large enough to absorb the last term in the right hand side \eqref{rhsw}, we obtain
\begin{eqnarray}
	&&\int_{\Omega}  e^{2s\varphi(0)}|fR(0)|^2 \,dx \nonumber\\
	&\leq& \dfrac 2{\sqrt s} \left(\int_{-T}^0\int_\Omega\left| P_1w\right|^2\,dxdt \right)^{\frac 12}
		\left(s\int_{-T}^0\int_ \Omega |\partial_t w|^2 \,dxdt\right)^{\frac 12} 
		+ M s^2\lambda^3 \int_{-T}^T\int_\Omega |w|^2\varphi^2 \,dxdt \nonumber\\
&\leq& \dfrac M{\sqrt s}\left( \int_{-T}^T\int_{\Omega} |Pw|^2 \,dxdt
+ s\int_{-T}^T\int_{\Gamma_0} \varphi |\partial_\nu w |^2d\sigma dt\right)\nonumber\\
&\leq& \dfrac M{\sqrt s}\left( \int_{-T}^T\int_{\Omega} e^{2s\varphi}|Lv|^2 \,dxdt
+ s \int_{-T}^T\int_{\Gamma_0}  e^{2s\varphi}|\partial_\nu v|^2d\sigma dt\right)\label{estim}
\end{eqnarray}
the weight function $\varphi$ being bounded from above and below.\\

From equation \eqref{WE4} and the properties \eqref{chi} of the cut-off function~$\chi$ and that of the weight $\varphi$ given in \eqref{psi} and \eqref{psi2}, one gets that
\begin{eqnarray*}
&&\int_{-T}^T\int_{\Omega}  e^{2s\varphi }|L v |^2\,dxdt \\
& \leq & M \int_{-T}^T\int_{\Omega}  e^{2s\varphi }|\chi f  \partial_t R |^2\,dxdt
+ M  \int_{-T}^T\int_{\Omega}  e^{2s\varphi }\left(|\chi' \partial_t z |^2 + |\chi'' z  |^2\right)\,dxdt\\
& \leq & M \int_{-T}^T\int_{\Omega}  e^{2s\varphi }|f |^2 |\partial_t R |^2\,dxdt
+ M  \left( \int_{-T}^{-T+\eta} + \int_{T-\eta}^T \right)\int_{\Omega}  e^{2s\varphi }\left(|\partial_t z|^2 + |z| ^2\right)\,dxdt\\
& \leq & M \int_{-T}^T\int_{\Omega}  e^{2s\varphi (0)}|f |^2 |\partial_t R |^2\,dxdt
+ M  \left( \int_{-T}^{-T+\eta} + \int_{T-\eta}^T \right)\int_{\Omega}  e^{2se^{\lambda C_0}}\left(|\partial_t z|^2 + |z| ^2\right)\,dxdt\\
& \leq & M \|\partial_t R\|_{L^2(-T,T;L^\infty(\Omega))} \int_{\Omega}  e^{2s\varphi (0)}|f |^2 \,dx\\
&&\qquad+~ M e^{2se^{\lambda C_0}} \left( \int_{-T}^{-T+\eta} + \int_{T-\eta}^T \right)\int_{\Omega}  \left(|\partial_t z|^2 + |z| ^2\right)\,dxdt.
\end{eqnarray*}
Using the energy estimate given in \eqref{energy}, and again the property \eqref{psi} of the weight $\varphi$, one gets 
\begin{eqnarray*}
\int_{-T}^T\int_{\Omega} e^{2s\varphi}|Lv|^2\,dxdt 
& \leq & M  \int_{\Omega} e^{2s\varphi(0)}|f|^2 \,dx
+ M e^{2se^{\lambda C_0}} \left( \int_{-T}^{-T+\eta} + \int_{T-\eta}^T \right) E_z(t) dt\\
& \leq & M   \int_{\Omega} e^{2s\varphi(0)}|f|^2 \,dx\\
&&+ M \eta  \left( ||R(0)||^2_{L^{\infty}(\Omega)} + ||R||^2_{H^1(0,T;L^{\infty}(\Omega))}\right) e^{2se^{\lambda C_0}}  \int_{\Omega} |f|^2 \,dx\\
& \leq &  M   \int_{\Omega} e^{2s\varphi(0)}|f|^2 \,dx.
\end{eqnarray*}
Gathering this last estimate with \eqref{estim}, we have proved
$$ \int_{\Omega}  e^{2s\varphi(0)}|fR(0)|^2 \,dx
 \leq \dfrac M{\sqrt{s}} \int_{\Omega} e^{2s\varphi(0)}|f|^2 \,dx
+ M\sqrt{s} \int_{-T}^T\int_{\Gamma_0} e^{2s\varphi}|\partial_\nu v|^2\,d\sigma dt.$$
Therefore, the assumption \eqref{RR} made on $R$ allow to obtain
$$ \int_{\Omega}  e^{2s\varphi(0)}|f|^2 \,dx\\
 \leq \dfrac M{\sqrt{s}}  \int_{\Omega}e^{2s\varphi(0)}|f|^2 \,dx
+ M\sqrt{s} \int_{-T}^T\int_{\Gamma_0} e^{2s\varphi}|\partial_\nu v|^2\,d\sigma dt,$$
and the choice of $s$ large enough gives
\begin{eqnarray*}
\int_{\Omega}  e^{2s\varphi(0)}|f|^2 \,dx
 &\leq &M\sqrt{s}  \int_{-T}^T\int_{\Gamma_0} e^{2s\varphi}|\partial_\nu v|^2\,d\sigma dt\\
 &\leq &M\sqrt{s}  \int_{-T}^T\int_{\Gamma_0} e^{2s\varphi}|\partial_\nu(\partial_t u)|^2\,d\sigma dt.
 \end{eqnarray*}
The proof of Theorem~\ref{Thm2} is then complete.
\endproof


\noindent{\bf Acknowledgements.} The author acknowledges Jean-Pierre Puel for this work accomplished together under is supervision during the masters degree in 2001. It has a been an important step toward the study of an inverse problem for the Schr\"odinger equation in the author's PhD and was also more recently strongly used in several works.\\

\bibliographystyle{plain}

\end{document}